\documentclass[12pt,a4paper]{article}
\usepackage{natbib}
\bibpunct{(}{)}{;}{a}{,}{,}
\usepackage{times}

\usepackage[english]{babel} %on pourrait ausi utiliser [english] \usepackage{latexSym}
\usepackage{amssymb}
\usepackage{amsmath}
\usepackage{amsfonts}
\usepackage[dvips]{graphicx}
\usepackage{subfigure}
\usepackage{epsfig}
\usepackage{overpic}
\usepackage{rotating}
\usepackage{here}
\usepackage{txfonts}
\usepackage{color}
\usepackage{tikz}

\newtheorem{theorem}{Theorem}[section]

\newtheorem{proposition}{Proposition}[section]

\newtheorem{definition}{Definition}[section]
\newtheorem{example}{Example}[section]
\newtheorem{remark}{\bf Remark}[section]

\newcommand{\bit}{\begin{itemize}}
\newcommand{\eit}{\end{itemize}}

\def\Prodi{\mathop{{\lower 3pt\hbox{\epsfxsize=15pt\epsfbox{pi.ps}}}}}
\def\prodi{\mathop{{\lower 1pt\hbox{\epsfxsize=8pt\epsfbox{pi.ps}}}}}
\input epsf.sty
\newcommand{\CAL}[1]{\mathcal{#1}}

\newcommand{\beq}{\begin{equation}}
\newcommand{\eeq}{\end{equation}}

\newcommand{\bea}{\begin{eqnarray}}
\newcommand{\eea}{\end{eqnarray}}
\newcommand{\beas}{\begin{eqnarray*}}
\newcommand{\eeas}{\end{eqnarray*}}

\newcommand{\BB}[1]{\mathbb{#1}}
\newcommand{\noi}{\noindent}

\begin{document}

\bigskip
\begin{center}
\textbf{\large{
Maxiset point of view for signal detection in inverse problems}}
\\[4mm]
\textbf{
Florent Autin$^1$, Marianne Clausel$^2$, Jean-Marc Freyermuth$^3$,  Cl\'ement Marteau$^4$ }
\\[1ex]
\today
\end{center}
$^1$ $^{3}$  Aix-Marseille Universit\'e, CNRS, Centrale Marseille, I2M, Marseille, France.\\
$^2$ Universit\'e de Lorraine, CNRS, Inria, IECL, F-54000 Nancy, France,\\
$^4$ Universit\'e Claude Bernard Lyon 1, CNRS UMR 5208, Institut Camille Jordan, F-69622 Villeurbanne, France.

%\tableofcontents

\section{Introduction}
  
\noi Over the last 20 years, the assessment of the performance of nonparametric function estimation methods mainly relied on the asymptotic minimax and oracle approaches.  More marginally used, the maxiset paradigm has been proved to be very useful to accurately describe the behaviour of some estimation procedures. In some cases, it allows to distinguish methods having comparable minimax performance. The question of adapting the maxiset concepts to the signal detection framework was often raised. This is the aim of this paper to rigorously extend this point of view to the signal detection framework and to discuss new related outcomes. \\

To this end, we will deal all along the paper with the Gaussian sequence space model 
\begin{equation}
y_k = b_k \theta_k + \varepsilon \xi_k, \quad k\in J,
\label{eq:model}
\end{equation}
where $(y_k)_{k\in J}$ denotes the observations, $J\subset \mathbb{N}^{*}$ is a subset of $\mathbb{N}^{*}$, $\theta = (\theta_k)_{k\in J}$ a non negative unknown sequence of interest, $(b_k)_{k\in J}$ a given sequence of non negative real numbers, $\varepsilon$ a noise level in $(0,1)$ and $(\xi_k)_{k\in J}$ a sequence of i.i.d. standard Gaussian random variables. The model (\ref{eq:model}) allows to describe several situations, as e.g. nonparametric regression or estimation of a function blurred by white noise. For more details on these models and their connection with the Gaussian sequence space model we refer the interested reader to \cite{Tsybakov}. For the sake of convenience, we will consider hereafter that $J= \mathbb{N}^*=\mathbb{N}\setminus\{0\}$.  We also stress that the model (\ref{eq:model}) allows to deal with so-called inverse problem models as described in \cite{Cavalier}. In such a setting, one is interested in doing inference on a function $f\in H$ in some Hilbert space $H$ from indirect and blurred observation of the form 
\begin{equation}
Y = Af + \varepsilon \xi,
\label{eq:ipmodel}
\end{equation}
where $A: H \rightarrow K$ denotes a compact operator acting from $H$ to another Hilbert space $K$, $\varepsilon$ a noise level and $\xi$ a Gaussian white noise. In particular, the sequence $(b_k^2)_{k\in \mathbb{N}^{*}}$ can be identified as the sequence of eigenvalues of the operator $A^\star A$ and the sequence $(\theta_k)_{k\in \mathbb{N}^{*}}$ as the one of the coefficients of $f$ in the singular values decomposition (SVD) basis associated to the operator $A$.   \\

For inverse problem models, the minimax paradigm has been widely used in order to assess the performance of estimation procedures. Roughly speaking, given a structural constraint on the vector $\theta$ of interest, typically of the form $\theta \in \Theta$ for some $\Theta \subset l_2(\mathbb{N}^{*})$, one measures the performance of a given estimator $\hat \theta$ through its maximal risk
$$ R(\hat\theta) = \sup_{\theta \in \Theta} \mathbb{E}_\theta \ \ell(\theta,\hat \theta),$$
where $\ell(.,.)$ denotes a given loss function. This paradigm has been widely used and discussed over the years. In several situations, a precise bound can be obtained on $R(\hat\theta)$ which allows to characterize how the maximal risk decreases with respect to the noise level $\varepsilon$.
 
More precisely, one can often exhibit a non decreasing positive sequence $(r_\varepsilon)_{\varepsilon>0}$, the so--called rate of convergence associated to the estimation procedure $\hat \theta$ with noise level $\epsilon$, such that $R(\hat\theta) \leq C r_{\varepsilon}$ for some positive constant $C$.  If this rate appears to be the smallest possible one, namely if there exists a positive constant $c$ such that 
$$ \inf_{\tilde \theta} \sup_{\theta \in \Theta} \mathbb{E}_\theta \ \ell(\theta,\tilde \theta) \geq c r_\varepsilon,$$
the sequence $(r_\varepsilon)_{\varepsilon>0}$ is called minimax rate of convergence over $\Theta$.  In the previous inequality, the infimum is taken over all possible estimators $\tilde \theta$ of $\theta$. We refer, e.g., to \cite{Tsybakov}, \cite{Jan11} for a non-exhaustive reference list.    \\

Under the minimax estimation paradigm, the performance of two given procedures can be compared through their respective rates of convergence according to a chosen functional set $\Theta$. However, it does not always allow for comparison if both procedures are 'minimax-optimal'. In addition the used criterium is quite pessimistic: the risk is measured at the slowest possible estimation precision over the set $\Theta$. Hence, it does not provide a fair comparison. To tackle these issues, an alternative point of view has been proposed in the seminal paper of \cite{KP2002}. The main idea can be stated as follows; given an estimation procedure $\hat \theta$ and a sequence of rates $(r_\varepsilon)_{\varepsilon>0}$, can we determine the set $\Theta_{MS}(\hat\theta)$ of sequences $\theta$ than are estimated by $\hat \theta$ at the rate  $(r_\varepsilon)_{\varepsilon>0}$? If yes, the set $\Theta_{MS}(\hat\theta)$ is called the maxiset associated to the procedure $\hat\theta$ for the rate  $(r_\varepsilon)_{\varepsilon>0}$. Under this paradigm, the best performing procedure, i.e., the 'maxiset-optimal' procedure, is the one whose associated maxiset strictly includes the maxisets of the others. Note that a very usual criticism concerns the situation where estimation methods have non nested maxisets. \cite{AFvS} discuss this important aspect of the maxiset approach explaining that, first, it is is somehow normal to find that some estimation methods are better in estimating some specific functions. In such a case, examining the 'form' of the maxiset will bring interesting information. Second, it may be possible to combine these procedures such that the maxiset combined procedure contains the union of the maxisets. The maxiset point of view has been generalized to various settings, see, e.g., \cite{A06}, \cite{RT08} or \cite{TW17}.   \\

\noi In the framework of signal detection, the minimax point of view has been widely investigated and was very fruitfully applied. We refer to \cite{Ingster_book}, \cite{Baraud}, \cite{ISS2012} or \cite{LLM2012} among others. Nevertheless, as in the estimation case, the minimax paradigm does not allow for a fully satisfying comparison between different testing procedures. The extension of the maxiset theory to this setting is a doorway to novel informative and rigorous math-stat study of these procedures. A flavor of the maxiset approach in signal detection framework has been discussed in \cite{ACF}. Nevertheless, the proper adaptation of this approach to the signal detection framework is a challenging problem that we tackle in this paper. We further discuss some new issues related to this theory. In particular, in this work, we aim at 
\begin{itemize}
\item highlighting the link between the space $\Theta$  and the sequence $(r_\varepsilon)_{\varepsilon>0}$ that both appear in the alternative hypothesis of the testing problem for different procedures based on $\chi^2$-statistics, 

\item comparing inverse and direct approaches in the light of the maxiset point of view.
%\item introduce and discuss the principle of maxiclass of operators, which will allow to characterize the kind of inverse problem that can be handled for a given situation.
\end{itemize}

In Section \ref{s:2}, we recall the minimax paradigm and then present the maxiset point of view for signal detection problems. Thereafter, in Section~\ref{s:3}, we state maxiset results  for both the inverse approach and the direct one (see Theorems $3.1$ and $3.2$). A crucial and perhaps surprising aspect in signal detection in inverse problem has been raised in \cite{LLM2011} where they compared direct and indirect testing procedures, ie., from the minimax point of view.  In Section~\ref{s:4} we succeed in comparing inverse and direct approaches in the light of the maxiset approach in many cases (see Proposition $4.1$) as for instance the moderately ill-posed inverse problem (see Proposition $4.2$). 
%In Section~\ref{s:3},  we compare  we focus on two important on $\chi^2$-based test statistics. In each case, we explicit the associated maxisets (see Theorems~\ref{maintheo1} and \ref{maintheo}), and provide a comparison of theses two procedures through this formalism. 

%In Section~ \ref{s:4}, we define in a formal way the notion of maxiclass of operator introduced above, and we make explicit the nature of the maxiclass for some specific inverse problems. 
Following a brief conclusion on the novelty of our results in Section~\ref{s:5}, we postpone in Section~\ref{s:6} all the related proofs.

\section{Signal detection in inverse problems}\label{s:2}

\subsection{The minimax paradigm for inverse problems in signal detection}

We consider the sequence space model (\ref{eq:model})
\[
y_k = b_k \theta_k +\varepsilon\xi_k, \quad k\in \mathbb{N}^{*}.
\]
The signal detection problem aims at determining whether or not the observations $y=(y_k)_{k\in \mathbb{N}^{*}}$ contains some signal. This question can be formalized as the following hypothesis testing problem: 
\begin{equation}
H_0:\theta = 0_{l_2(\mathbb{N}^*)} \quad \mathrm{against} \quad H_1:\theta \in \Theta, \: \left\|\theta\right\|^2 \geq  r_\varepsilon^2,
\label{eq:testingpb}
\end{equation}
for some non decreasing positive sequence $r=(r_\varepsilon)_{\varepsilon>0}$ depending on $\varepsilon$.\\
\noi In the alternative hypothesis $H_1$, the set $\Theta$ denotes a subset of $l_2(\mathbb{N}^{*})$. The requirement $\theta \in \Theta$ can be thought either as a structural constraint on the signal or as a regularity condition on the underlying function $f$ in model (\ref{eq:ipmodel}). In the same time, the constraint $\left\|\theta\right\|^2 \geq  r_\varepsilon^2$ corresponds to an energy condition that allows to quantify the amount of signal available in the observations. Another problem closely related to signal detection is pattern recognition. In this
case, one aims at testing the adequation between the observations and a given reference signal $\theta^{o}$. Having in mind that up to the change of variable $\theta\leftarrow \theta-\theta^{o}$, these two problems are equivalent we shall only focus in the sequel on the signal detection problem.\\

In the sequel we denote as $\Delta = \Delta(y)$ a testing procedure. It is a measurable function of $y=(y_k)_{k\in \mathbb{N}^{*}}$, such that $\Delta \in \lbrace 0,1 \rbrace$, with the convention that we reject $H_0$ if $\Delta = 1$ - that corresponds to accept $H_1$ - and do not reject the null hypothesis $H_0$ otherwise. Given a level $\alpha \in (0,1)$, $\Delta = \Delta_\alpha$ is called a level-$\alpha$ test if and only if 
$$ \mathbb{P}_{0_{l_2(\mathbb{N}^*)}}(\Delta_\alpha = 1) \leq \alpha.$$
The risk under the alternative hypothesis is often measured through the maximal Type-II error over the set $\Theta$, as 
$$\beta(\Delta_\alpha,\Theta,r_\varepsilon):=  \sup_{\theta \in \Theta, \:  \left\| \theta\right\|^2 \geq r_\varepsilon^2   }\mathbb{P}_{\theta}(\Delta_{\alpha} = 0).$$
In particular, given a level $\beta \in (0,1)$, the level-$\alpha$ test $\Delta_\alpha$ is said to be powerful if its maximal Type-II error can be bounded by $\beta$, namely $\beta(\Delta_\alpha,\Theta,r_\varepsilon) \leq \beta$. In this context, the minimax paradigm has been at the core of several investigations over the last decades. Given both $\alpha$ and $\beta \in (0,1)$ some fixed levels and a given set $\Theta$, the separation rate $R_\varepsilon(\Delta_\alpha,\Theta)$ associated to a given level $\alpha$-test $\Delta_\alpha$ is defined as 
$$ R_\varepsilon(\Delta_\alpha,\Theta) = \inf \left\lbrace r_\varepsilon >0, \ \beta(\Delta_\alpha,\Theta,r_\varepsilon) \leq \beta \right\rbrace.$$
Typically, we expect that $R_\varepsilon(\Delta_\alpha,\Theta) \rightarrow 0$ as $\varepsilon \rightarrow 0$, although it strongly depends on the considered setting. The minimax separation rate associated to the testing problem (\ref{eq:testingpb}) for a given set $\Theta$ is then defined as $r^\star = (r^\star_\varepsilon)_{\varepsilon>0}$ where
$$ r_\varepsilon^\star := \inf_{\Delta} R_\varepsilon(\Delta,\Theta),$$
where the infimum is taken over all possible level-$\alpha$ tests $\Delta$. We refer to \cite{Ingster_book} or \cite{Baraud} for exhaustive discussions on theses definitions. Determining the minimax rate of convergence for a given problem is quite informative. For substantial account on the subject, see e.g., \cite{Ingster_book}, \cite{Baraud}, \cite{Butucea}, \cite{LLM2012} or \cite{Lacour2014} among others.

\subsection{Inverse and direct approaches}

%Following the context, and in particular the structure of the set $\Theta$, several testing procedures have been proposed in order to deal with the testing problem (\ref{eq:testingpb}). Our aim in this section is not to provide an exhaustive list. We will focus on two particular tests based on $\chi^2$ statistics that have been proved to provide interesting performance in some minimax settings. \\
Several testing procedures have been proposed in order to deal with the testing problem (\ref{eq:testingpb}). In this section we will focus on two $\chi^2$-based test statistics that have been proved to perform well in some minimax settings. \\

The amount of signal contained in the observations can be measured through the quantity $\| \theta \|^2$. Then, a natural way to provide a decision rule is to construct an estimator of this quantity. \\

\underline{Inverse approach} (IP): According to the sequence model (\ref{eq:model}), each coefficient $\theta_k$ can be estimated by $b_k^{-1} y_k$ provided that the sequence of $b=(b_k)_{k \in \mathbb{N}^*}$ is known. For a given non decreasing sequence of non negative integers $(D_\varepsilon)_{\varepsilon>0}$, this leads to the testing procedure  $\Delta_{\alpha,\varepsilon}^{IP}$
\begin{equation}
\Delta_{\alpha,\varepsilon}^{IP} = \mathbf{1}_{\lbrace T_{D_\varepsilon} > t_{\alpha,\varepsilon} \rbrace},\text{ where } T_{D_\varepsilon}= \sum_{k=1}^{D_\varepsilon} b_k^{-2}(y_k^2 - \varepsilon^2)
\label{eq:inverse_test}
\end{equation}

and $t_{\alpha,\varepsilon}$ is a threshold value that allows to control the Type-I error of the test. The integer $D_\varepsilon$ plays a similar role than a regularization parameter in estimation. According to the choice of the set $\Theta$, specific choices for $D_\varepsilon$ are available. We refer, e.g., to \cite{LLM2012} for more details. We stress that a weighted variant of this procedure has been proposed and investigated in several papers, as e.g. \cite{ISS2012}, allowing to obtain sharp asymptotic results. \\

\underline{Direct approach} (DP): Since the sequence model (\ref{eq:model}) is derived from the model (\ref{eq:ipmodel}), we remark that the test (\ref{eq:inverse_test}) is essentially based on an inversion of the operator $A$ at hand. Such an inversion appears to be quite natural in an estimation context in which we provide a reconstruction of the unknown function $f$. In a signal detection framework, this is no longer required. Indeed, setting $\vartheta = b.\theta$ (i.e. $\vartheta_k = b_k \theta_k$ for all $k \in \mathbb{N}^{*}$), we can remark that both assertion $\theta = 0_{l_2(\mathbb{N}^*)}$ and $\vartheta =0_{l_2(\mathbb{N}^*)}$ are equivalent. In other words, the testing problem, for some non decreasing positive sequence $\mu=(\mu_\varepsilon)_{\varepsilon> 0}$:
\begin{equation}
\tilde{H}_0: \vartheta=0_{l_2(\mathbb{N}^*)} \quad \mathrm{against} \quad \tilde{H}_1 : \vartheta \in \tilde{\Theta}, \ \| \vartheta \|^2 \geq \mu_{\varepsilon}^2,
\label{eq:testing_direct}
\end{equation}
for some set $\tilde{\Theta}$ only differs from (\ref{eq:testingpb}) by the alternative. In some sense, (\ref{eq:testing_direct}) does not take into account the fact that the data are distorted by a compact operator: we threat the data as a `direct' problem and deal with a model of the form
\begin{equation}
Y = g + \varepsilon \xi \quad \mathrm{or} \quad y_k = \vartheta_k + \varepsilon \xi_k \quad \forall k\in \mathbb{N}^{*},
\label{eq:modeldir}
\end{equation}
where $g = Af \in K$.  Consequently, for a given non decreasing sequence of non negative integers $(D_\varepsilon)_{\varepsilon>0}$ we can introduce the test $\Delta_{\alpha,\varepsilon}^{DP}$ as 
\begin{equation}
\Delta_{\alpha,\varepsilon}^{DP} = \mathbf{1}_{\lbrace S_{D_\varepsilon} > s_{\alpha,\varepsilon}\rbrace} \quad \mathrm{where} \quad %S_{D_\varepsilon} = \sum_{k=1}^{D_\varepsilon} y_k^2,
S_{D_\varepsilon} =\sum_{k=1}^{D_\varepsilon} (y_k^2 - \varepsilon^2),
\label{eq:direct_test}
\end{equation}
and $S_{D_\varepsilon}$ corresponds to an estimator of $\| \vartheta \|^2$ and $s_{\alpha,\varepsilon}$ denotes an appropriate threshold, allowing a control of the Type-I error.
This test provides interesting performance when dealing with (\ref{eq:testing_direct}). But surprisingly, this is also the case for the testing problem (\ref{eq:testingpb}) in some specific situations. We refer to \cite{LLM2011} for an extensive discussion on the subject with a minimax point of view. One of the aim of this paper is to complete this discussion using a maxiset point of view. This notion is extended to the signal detection context in the next section.

\subsection{The maxiset point of view in signal detection problem}

In nonparametric function estimation, an estimator is a sequence $(\hat \theta_{D_\varepsilon})_{\varepsilon>0}$ possibly indexed by a regularization (smoothing) parameter $D$ and by the noise level $\varepsilon$, is a sequence $(\hat \theta_{D_\varepsilon})_{\varepsilon>0}$.
In the minimax setting, given a functional set $\Theta$, we determine sequence of rates $r=(r_\varepsilon)_{\varepsilon>0}$ such that for any $\varepsilon \in (0,1),$
\begin{equation}
\inf_{\hat{\theta}}\sup_{\theta \in \Theta} \BB{E}_\theta \|\hat{\theta} - \theta \|^2 \geq c r_\varepsilon^2 \text{ and } \sup_{\theta \in \Theta} \BB{E}_\theta \|\hat{\theta}_{D_\varepsilon} - \theta \|^2 \leq C r_\varepsilon^2,
\label{eq:control}
\end{equation}
for some constant $c, C>0$. Under the maxiset paradigm, we are giving a sequence of rates $r=(r_{\varepsilon})_{\varepsilon>0}$, and we exhibit the largest functional set $\Theta \subset l_2(\mathbb{N}^{*})$ for which (\ref{eq:control}) holds. \\

We will now adapt the maxiset point of view to the signal detection framework. Given a sequence $r=(r_\varepsilon)_{\varepsilon>0}$, we determine the largest set for which the maximal Type-II error can be controlled with a prescribed error of our testing problem. This is formalized in the following definition. 

\begin{definition}
\label{def:ms}
For a fixed $(\alpha, \beta)\in (0,1/2)^2$, let $\Delta_{\alpha}=(\Delta_{\alpha,\varepsilon})_{\varepsilon>0}$ be a sequence of testing procedures and $(r_\varepsilon)_{\varepsilon>0}$ a decreasing sequence of non negative real numbers. The maxiset of $\Delta_{\alpha}$ associated to the separation rate $r=(r_\varepsilon)_{\varepsilon>0}$, is the largest sequence space $\Theta$ in $l_2(\mathbb{N}^{*})$  such that, for all $\varepsilon \in (0,1)$
\beas
\mathbb{P}_{0_{l_2(\mathbb{N}^*)}}(\Delta_{\alpha,\varepsilon} = 1)\leq \alpha \: \text{ and } \sup_{\theta \in \Theta, \: \left\|\theta\right\|^2\geq r_\varepsilon^2}\mathbb{P}_\theta(\Delta_{\alpha,\varepsilon}=0)\leq \beta.
\eeas
\end{definition}

\noi This definition can be generalized in a straightforward way to the testing problem (\ref{eq:testing_direct}). In the following, we denote the maxiset as $MS(\Delta_{\alpha},r):=MS( (\Delta_{\alpha,\varepsilon})_{\varepsilon>0}, (r_\varepsilon)_{\varepsilon>0})$. Note that it clearly corresponds to the following set:
\bea \label{maxidef-r}
MS( \Delta_{\alpha}, r) =\left\{\left(\theta_k\right)_{k\in \mathbb{N}^{*}} \in l_2(\mathbb{N}^{*}) :\forall \varepsilon \in (0,1),   \left[\left\|\theta\right\|^2\geq r_\varepsilon^2 \Rightarrow \mathbb{P}_\theta\left[\Delta_{\alpha, \varepsilon} = 0 \right]\leq \beta \right]\right\}.
\eea

In Section~\ref{s:3}, we shall derive an explicit expression of the maxisets for some tests based on $\chi^2$ statistics (see Theorem~\ref{maintheo1} and \ref{maintheo2} below). 

Following Definition \ref{def:ms} we let the reader be convinced that, for a given sequence of testing procedures $(\Delta_{\alpha, \varepsilon})_{\varepsilon>0}$, there is an embedding result between its maxisets associated with different choices of detection rates. 

\begin{proposition} \label{pro:emb}
Let $(r'_\varepsilon)$ and $(r_\varepsilon)$ two sequences of detection rates such that $r_\varepsilon=o(r'_\varepsilon)$ as $\varepsilon\to 0$. Consider a sequence of testing procedures $\Delta_{\alpha}=(\Delta_{\alpha, \varepsilon})_{\varepsilon>0}$. Then, for some $C>0$
\[
MS(\Delta_{\alpha}, r) \subset MS(\Delta_{\alpha}, Cr'). 
\]
\end{proposition}

The previous embedding entails that the set of detectable functions increases as soon as we relax the constraint on the lowest minimal energy required. 

%We now give some details about the test procedures that we shall consider as well as the explicit expression of their associated maxisets.

\section{Testing procedures and their maxiset performance}\label{s:3}

In this section, we first provide a description of the maxisets associated respectively to the procedures $(\Delta_{\alpha,\varepsilon}^{IP})_{\varepsilon>0}$ and $(\Delta_{\alpha,\varepsilon}^{DP})_{\varepsilon>0}$ defined in ~\eqref{eq:inverse_test} and \eqref{eq:direct_test}. From now on, we fix the two thresholds $t_{\alpha,\varepsilon}$ and $s_{\alpha,\varepsilon}$ involved in the definition of these testing procedures as
\begin{equation}
    t_{\alpha,\varepsilon} = C_{\alpha,1}\, \varepsilon^2\sqrt{\sum_{k=1}^{D_\varepsilon}b_k^{-4}} \quad \mathrm{and} \quad  s_{\alpha,\varepsilon}= C_{\alpha,2}\, \varepsilon^2 \sqrt{D_\varepsilon} \quad \forall \varepsilon \in (0,1),
    \label{eq:hyp-t}
\end{equation}
where $C_{\alpha,1}, C_{\alpha,2}>0$ denote two  explicit constants that guarantee that the considered procedures have a Type-I error controlled by $\alpha$ for all $\varepsilon \in (0,1)$. Undoubtedly,  the smaller $\alpha$ is the bigger $C_{\alpha,1}$ and $C_{\alpha,2}$ are. For the sake of simplicity, we do not use the $1-\alpha$ quantile of the respective test statistics $(S_{D_\varepsilon})_{\varepsilon>0}$ and $(T_{D_\varepsilon})_{\varepsilon>0}$. Such a change will not modify the spirit of the results displayed in this paper but will induce more technical details.  \\
%and that
%\begin{equation}
%s_{\alpha,\varepsilon}\sim \varepsilon^2 \sqrt{D_\varepsilon}.%2 \varepsilon^4 \sqrt{D_\varepsilon}.
%  \label{eq:hyp-s}
%\end{equation}

Below we characterize the maxisets associated to the considered procedure for general separation rates. In such a setting, these sets are poorly informative but they highlight valuable information on the problem provided that we impose some structural constraints.

\subsection{A general characterization of the maxisets}\label{s:3-maxi}

We start our investigations with a general description of the maxisets associated to the procedures $(\Delta_{\alpha,\varepsilon}^{IP})_{\varepsilon>0}$ and $(\Delta_{\alpha,\varepsilon}^{DP})_{\varepsilon>0}$ for any chosen rate of detection $(r_\varepsilon)_{\varepsilon>0}$. To this end, we introduce the two sets $\CAL{F}_{r, D}(C)$ and $\CAL{G}_{r, D}(C) $ which will be of first importance in the sequel.

\begin{definition}%[sequence spaces]
Let $(D_\varepsilon)_{\varepsilon>0}$ be an increasing sequence of non negative integers. Let $(r_\varepsilon)_{\varepsilon>0}$ be a decreasing sequence of non negative real numbers . For any $C>0$ we set
\beas
\CAL{F}_{r,D}(C) = \left\{\theta\in l_2(\mathbb{N}^{*}), \:   \forall \varepsilon\in (0,1) ; \left[\|\theta\|^2\geq r_\varepsilon^2\Rightarrow \sum_{k=1}^{D_\varepsilon}\theta_k^2 > C\varepsilon^2\sqrt{\sum_{k=1}^{D_\varepsilon}b_k^{-4}}\right] \right\},
\eeas

\beas
\CAL{G}_{r,D}(C) = \left\{\theta\in l_2(\mathbb{N}^{*}), \:     \forall \varepsilon\in (0,1) ; \left[\|\theta\|^2\geq r_\varepsilon^2\Rightarrow \sum_{k=1}^{D_\varepsilon}b_k^2\theta_k^2>  C\varepsilon^2\sqrt{D_\varepsilon}\right] \right\}.
\eeas
\end{definition}

The following results emphasizes that these sets provide a characterization of the maxisets associated to the procedures $(\Delta_{\alpha,\varepsilon}^{IP})_{\varepsilon>0}$ and $(\Delta_{\alpha,\varepsilon}^{DP})_{\varepsilon>0}$. 

\begin{theorem}\label{maintheo1}
Consider $(\alpha,\beta)\in (0,1/2)^2$. Let $(t_{\alpha,\varepsilon})_{\varepsilon>0}$ and $(s_{\alpha,\varepsilon})_{\varepsilon>0}$ satisfying  \eqref{eq:hyp-t}. Consider the two sequences of testing procedures $(\Delta^{IP}_{\alpha,\varepsilon})_{\varepsilon>0}$ and $(\Delta^{DP}_{\alpha,\varepsilon})_{\varepsilon>0}$ defined respectively in \eqref{eq:inverse_test} and \eqref{eq:direct_test}. 
We have the two following maxiset results for any choice of detection rates $r=(r_{\varepsilon})_{\varepsilon>0}$ and $\mu=(\mu_{\varepsilon})_{\varepsilon>0}$:  
\begin{enumerate}
\item There exist two positive constants $C_{\min}(\alpha,\beta)$ and $C_{\max}(\alpha,\beta)$ depending on $C_{\alpha,1}$ and $\beta$ such that:
$$\CAL{F}_{r,D}(C_{\max}(\alpha,\beta))\subset MS(\Delta_{\alpha}^{IP}, r) \subset \CAL{F}_{r,D}(C_{\min}(\alpha,\beta)),$$
rewritten as: $MS(\Delta_{\alpha}^{IP}, r) = \CAL{F}_{r,D}$.
\item There exist two positive constants $C'_{\min}(\alpha,\beta)$ and $C'_{\max}(\alpha,\beta)$ depending on $C_{\alpha,2}$ and $\beta$ such that:
$$\CAL{G}_{\mu,D}(C'_{\max}(\alpha,\beta))\subset MS(\Delta_{\alpha}^{DP}, \mu) \subset \CAL{G}_{\mu,D}(C'_{\min}(\alpha,\beta)),$$
rewritten as: $MS(\Delta_{\alpha}^{DP}, \mu) = \CAL{G}_{\mu,D}$.
\end{enumerate}
\end{theorem}

\begin{remark} In Section $6$, we provide explicit values of the constants $C_{\max}(\alpha,\beta)$ and $C_{\min}(\alpha,\beta)$ (see (\ref{Cmax}) and (\ref{Cmin})). The constants $C'_{\max}(\alpha,\beta)$ and $C'_{\min}(\alpha,\beta)$ are obtained from the values of $C_{\max}(\alpha,\beta)$ and $C_{\min}(\alpha,\beta)$ by replacing $C_{\alpha,1}$ by $C_{\alpha,2}$.
\end{remark}
Surprisingly, the maxisets in the testing case have a completely different form compared to results obtained in the estimation case. Indeed, according to \cite{KP2002}, the constraint that a given procedure attains the rate $(r_\varepsilon)_{\varepsilon>0}$ induces a tail constraint on the signal of interest in the estimation problem. This is no more the case in the signal detection problem. The Theorem \ref{maintheo1} above indicates that the procedures $(\Delta^{IP}_{\alpha,\varepsilon})_{\varepsilon>0}$ and $(\Delta^{DP}_{\alpha,\varepsilon})_{\varepsilon>0}$ are able to detect only signals satisfying, for any $\varepsilon\in (0,1)$ such that $\|\theta\|^2 \geq r_\varepsilon^2 $
\begin{equation} 
\sum_{k=1}^{D_\varepsilon} \theta_k^2 > C \varepsilon^2 \sqrt{\sum_{k=1}^{D_\varepsilon} b_k^{-4}}, \quad \forall \varepsilon\in (0,1), 
\label{eq:contrainte}
\end{equation}
\begin{equation} 
\sum_{k=1}^{D_\varepsilon} b_k^2\theta_k^2 > C \varepsilon^2 \sqrt{D_\varepsilon}, \quad \forall \varepsilon\in (0,1),
\label{eq:contrainte}
\end{equation}
and for $C$ large enough. In particular, the constraint (\ref{eq:contrainte}) indicates that there should be enough signal on the frequencies investigated by the test statistics $T_{D_\varepsilon}$. Nothing is said regarding the high frequencies, i.e. the coefficient after the rank $D_\varepsilon$. The maxiset results of Theorem \ref{maintheo1}, constrasts with usual ones, since one does not describe maxisets in terms of smoothness spaces. Moreover, this constraint has already been highlighted in, e.g. \cite{LLM2012}. Hence, the maxiset paradigm is poorly informative in such a context. In the following section, we prove that an additional structural assumption on the maxiset provides valuable informations on the signal that can be detected by the procedures we are interested in.

\subsection{A robust version of maxisets for tests}
\label{s:robust}

The main spirit of the previous section is that the functions that can be detected by the procedures $(\Delta^{IP}_{\alpha,\varepsilon})_{\varepsilon>0}$ and $(\Delta^{DP}_{\alpha,\varepsilon})_{\varepsilon>0}$ have enough energies for low frequencies. It means in particular that our testing procedures are very sensitive to the trend of the signal. In what follows, we shall require some robustness of our procedure with respect to this low frequency part of the signal, provided we have enough information. Indeed, in many practical situation, the signal is preprocessed or filtered and we want to have theoretical guarantees about signal detection remaining still valid in this context. 

This structural constraint on the maxiset of interest, can be reformulated in a more formal way as follows:  

\begin{definition}
A set $\mathcal{H} \subset l_2(\mathbb{N}^{*})$ satisfies the decimation constraint if 
$$ \theta \in \mathcal{H} \quad \Rightarrow \quad \theta^{(-n)} \in \mathcal{H} \ \  \forall n\in \mathbb{N}^{*}, $$
where for any $\theta \in l_2(\mathbb{N}^{*})$ and $n\in \mathbb{N}^{*}$, $\theta^{(-n)} = (\theta_k \mathbf{1}_{\lbrace k> n \rbrace})_{k\in \mathbb{N}^{*}}$.
\end{definition}

%{\color{red} Justification heuristique de cette contrainte?? }
We stress that such a condition is for instance satisfied by all the sets $\Theta \subset l_2(\mathbb{N}^*)$ of the form 
$$ \Theta = \left\lbrace v \in l_2(\mathbb{N}^*): \sum_{k=1}^{+\infty} w_k v_k^\gamma < L \right\rbrace,$$
for some positive sequence $(w_k)_{k\in \mathbb{N}^*}$ and positive constants $\gamma,L$. %Such a set traduces some smoothness of interest through the decay of the coefficient of the target of interest. \\
Such a set describes some smoothness conditions through the decay of the coefficients of the function of interest.

Hereafter we define two sequence spaces that are restriction of $\CAL{F}_{D}$ and $\CAL{G}_{D}$ to sequences satisfying the decimation constraint and that depend on the chosen detection rates $r=\left(r_\varepsilon\right)_{\varepsilon>0}$ and $\mu =\left(\mu_\varepsilon\right)_{\varepsilon>0}$ appearing in the definition of maxiset given in (\ref{maxidef-r}).

\begin{definition}%[sequence spaces]
Let $r=\left(r_\varepsilon\right)_{\varepsilon>0}$ be a decreasing sequence of non negative real numbers and $D=(D_\varepsilon)_{\varepsilon>0}$ be an increasing sequence of non negative integers. For any $C>0$ we set
\beas
\CAL{F}^{dec}_{r, D}(C) = \left\{\theta\in l_2(\mathbb{N}^{*}), \: \forall \varepsilon\in (0,1) ; \sum_{k>D_\varepsilon}\theta_k^2 < r_{\varepsilon}^2 - C\varepsilon^2\sqrt{\sum_{k=1}^{D_\varepsilon}b_k^{-4}} \right\},
\eeas

\beas
\CAL{G}^{dec}_{r, D}(C) = \left\{\theta\in l_2(\mathbb{N}^{*}), \:   \forall \varepsilon\in (0,1) ; \sum_{k>D_\varepsilon}b_k^2\theta_k^2 < r_{\varepsilon}^2 -C\varepsilon^2\sqrt{D_\varepsilon} \right\}.
\eeas
\end{definition}

%If we constrain the maxiset to satisfy the decimation condition, 
If we are searching for the largest set, satisfying both the requirement of Definition~\ref{def:ms} and robust with respect to decimation, we define the so-called robust maxisets and we retrieve exactly the sets introduced in the previous definition, up to some constants. \\

Undoubtedly, following the definition of $D=(D_{\varepsilon})_{\varepsilon>0}$, $r=(r_{\varepsilon})_{\varepsilon >0}$ and $C$ the sequence spaces above can be identical to the empty space. In the sequel, we especially focus on the cases where these sequences spaces are not the empty space.

\begin{definition} \label{admi} For any chosen sequences $D=(D_{\varepsilon})_{\varepsilon>0}$, $r=(r_{\varepsilon})_{\varepsilon >0}$ and $C>0$, we say that $(r,D,C)$ is $\CAL{F}$-admissible (respectively $\CAL{G}$-admissible) if and only if $\CAL{F}^{dec}_{r, D}(C)$ (respectively $\CAL{G}^{dec}_{r, D}(C)$) is not the empty space.
\end{definition}

\begin{remark}
Following  Definition \ref{admi}, note that $(r,D,C)$ is $\CAL{F}$-admissible and $\CAL{G}$-admissible for rates of detection that do not converge to zero too fast as $\varepsilon$ tends to zero.
\end{remark}

\begin{theorem}\label{maintheo2}
Consider $(\alpha,\beta)\in (0,1/2)^2$. Let $(t_{\alpha,\varepsilon})_{\varepsilon>0}$ and $(s_{\alpha,\varepsilon})_{\varepsilon>0}$ satisfying \eqref{eq:hyp-t}. Consider the two sequences of testing procedures $(\Delta^{IP}_{\alpha,\varepsilon})_{\varepsilon>0}$ and $(\Delta^{DP}_{\alpha,\varepsilon})_{\varepsilon>0}$ defined respectively in \eqref{eq:inverse_test} and \eqref{eq:direct_test}. 
%Consider that the maxiset must satisfy the decimation constraint. 
We then respectively denote them by  $MS^{dec}(\Delta^{IP}_{\alpha}, r)$ and  $MS^{dec}(\Delta^{DP}_{\alpha}, \mu)$ the respective robust maxisets associated with the chosen rates $r=(r_{\varepsilon})_{\varepsilon>0}$ and $\mu=(\mu_{\varepsilon})_{\varepsilon>0}$. We have the following maxiset results:
\begin{enumerate}
\item If $(r,D,C_{\max}(\alpha,\beta))$ is $\CAL{F}$-admissible, then:
\begin{equation}\label{eq:robMS1a}
\CAL{F}^{dec}_{r,D}(C_{\max}(\alpha,\beta))\subset MS^{dec}(\Delta_{\alpha}^{IP}, r) \subset  \CAL{F}^{dec}_{\sqrt{2}r,D}(C_{\min}(\alpha,\beta)),
\end{equation}
rewritten as: $MS^{dec}(\Delta_{\alpha}^{IP}, r) = \CAL{F}_{r,D}^{dec}.$
\item If $(\mu,D,C'_{\max}(\alpha,\beta))$ is $\CAL{G}$-admissible, then:\begin{equation}\label{eq:robMS2a}
\CAL{G}^{dec}_{\mu,D}(C'_{\max}(\alpha,\beta))\subset MS^{dec}(\Delta_{\alpha}^{DP}, \mu)\subset \CAL{G}^{dec}_{\sqrt{2}\mu,D}(C'_{\min}(\alpha,\beta)),
\end{equation}
rewritten as: $MS^{dec}(\Delta_{\alpha}^{DP}, \mu) = \CAL{G}_{\mu,D}^{dec}.$
\end{enumerate}
\end{theorem}

\begin{remark}
The constants stated in Theorem \ref{maintheo2} are similar to those in Theorem \ref{maintheo1}.
\end{remark}
We observe that, as in the estimation case, the maxiset  with a decimation constraint depends on the tail of the sequence $(\theta_k)_{k \in \BB{N}^*}$ of interest. Note that in the framework of signal detection, the situation is much more intricate than in estimation since one has several parameters to deal with: the rate of convergence $(r_\varepsilon)_{\varepsilon>0}$, the nature of the operator eventually involved in the inverse signal problem detection and $\beta\in (0,1)$ the Type-II error that has to be controlled.

According to the relative growth of the levels of possible energies of the signal and the sums of the power of the eigenvalues of the operator involved in the signal detection problem, the nature of the maxiset related to the sequence of testing procedures might be different. Consider the case where the sequence of testing procedures is $\Delta_{\alpha,\varepsilon}^{IP}$. There are two extreme situations:
\begin{itemize}
%\item First case: $r_\varepsilon^2\ll\varepsilon^2\sqrt{\sum_{1}^{D_\varepsilon}b_k^{-4}}$ as $\varepsilon\to 0$. In this case, Theorem~\ref{maintheo2} implies that the maxiset is empty, which means that our testing procedures are unable to detect any signal whatever its energy is.
\item First case: $r_\varepsilon^2\ll\varepsilon^2\sqrt{\sum_{k=1}^{D_\varepsilon}b_k^{-4}}$ as $\varepsilon\to 0$. In this case, Theorem~\ref{maintheo2} implies that the robust maxiset is empty. It means that under the considered noise level $\varepsilon$, whatever the signal we consider is, our procedure is never robust under decimation.  %we cannot satisfy both the decimation condition and have good statistical performances for our testing procedure. It means that, over the required performance level, our testing procedures will be always sensitive to the acquisition process whatever the energy of the signal is. \textcolor{red}{FLOU : MESSAGE PAS CLAIR}

%\item Second case: $r_\varepsilon^2\gg\varepsilon^2\sqrt{\sum_{1}^{D_\varepsilon}b_k^{-4}}$ as $\varepsilon\to 0$. In this case, the maxiset is non empty and does not depend on the operator. In particular, provided that we have enough energy in our signal, the performance of our detection procedure does not depend on the underlying inverse problem we are considering.  

\item Second case: $r_\varepsilon^2\gg\varepsilon^2\sqrt{\sum_{k=1}^{D_\varepsilon}b_k^{-4}}$ as $\varepsilon\to 0$. In this case, the robust maxiset is non empty and does not depend on the operator. In particular, provided that we have enough energy in our signal, the performance of our detection procedure does not depend on the underlying inverse problem we are considering.  
\end{itemize}
The transition case where the two sequences  $r_\varepsilon^2$ and $\varepsilon^2\sqrt{\sum_{1}^{D_\varepsilon}b_k^{-4}}$ are equally balanced. Here, the robust maxiset can be explicitly embedded as follows
\[
MS^{dec}(\Delta_{\alpha}^{IP}, r) \subset\left\{\theta\in l_2(\mathbb{N}^{*}), \:  \forall \varepsilon\in (0,1) ;\sum_{k>D_\varepsilon}\theta_k^2 <  r_\varepsilon^2 \right\},
\]
for some positive constant $C$. The set on the right-hand side of the previous embedding provides a control on the tail of the sequence of interest by the considered rate $(r_\varepsilon)_{\epsilon>0}$.  In particular, the faster the sequence $(r_\varepsilon)_{\varepsilon>0}$ converges toward $0$, the smoother the detectable function. \\

Note that since the sequence space $\mathcal{F}^{dec}_{r,D}(C)$ may be empty or not, depending on the value of $C$, one cannot conclude that $MS(\Delta^{IP}_{\alpha}, r)$ is non-empty in whole generality. Similar comments are also valid if we considering the sequence of testing procedures $\Delta_{\alpha,\varepsilon}^{DP}$ and if we compare the behavior of the two sequences $(\mu_\varepsilon^2)$ and $(\varepsilon^2\sqrt{D_\varepsilon})$ in the sequence space $\mathcal{F}^{dec}_{r,D}(C)$ as $\varepsilon\to 0$.\\

\section{Comparison of direct and inverse approaches}\label{s:4}

In this section, we will take advantage of tools developed in the previous section in order to compare the \textit{direct} and \textit{inverse} approaches in a signal detection. \\

Indeed, we have seen that both problems (\ref{eq:testingpb}) and (\ref{eq:testing_direct}) only differ by there alternatives. The tests $(\Delta_{\alpha,\varepsilon}^{IP})_{\varepsilon>0}$ and $(\Delta_{\alpha,\varepsilon}^{DP})_{\varepsilon>0}$ have been specially designed in order to answer separately to each of these problems. Now, a challenging question is to compare the alternative and to check whether the inverse (resp. direct) approach is pertinent for the problem (\ref{eq:testing_direct}) (resp. (\ref{eq:testingpb})). This comparison will be provided under the maxiset paradigm, using the robust version displayed in Section \ref{s:robust} above. To improve readability of our results, we now denote $MS$ instead of $MS^{dec}$, when denoting the robust maxisets.\\

In order to provide a fair comparison between both testing procedure, we have to state a dependency between the rates $r=(r_\varepsilon)_{\varepsilon>0}$ and $\mu =(\mu_\varepsilon){\varepsilon>0}$. Indeed, both procedures do not come up in the same space. For instance, in the minimax paradigm, the rates are often faster for 'direct' alternative than for the inverse case.  Below, we fix this dependency according to previous calibration that has been investigated in the minimax paradigm (see, e.g. \cite{LLM2011}). Concerning the regularization parameter $(D_\varepsilon)_{\varepsilon>0}$, we will keep the same value for both testing procedure: the idea is to work with the same number of coefficients (same amount of information).

\begin{proposition}\label{maincor}
Fix $(\alpha, \beta) \in (0,1/2)^{2}$. Choose $r=(r_\varepsilon)_{\varepsilon>0}$, $\mu=(\mu_\varepsilon)_{\varepsilon>0}$ and $D=(D_\varepsilon)_{\varepsilon>0}$ such that, for any $\varepsilon  \in (0,1)$, $\mu_\varepsilon = b_{D_\varepsilon}\,r_\varepsilon$. Then, provided that 
\begin{equation} 
C'_{\max}(\alpha, \beta) \sqrt{k} \leq C_{\min}(\alpha, \beta) b_k^2 \sqrt{ \sum_{j=1}^k b_j^{-4}} \quad  \hbox{for all }  \: k \in \mathbb{N}^*, 
\label{eq:cond}
\end{equation}
we get
\bea\label{eq:embIPDP}
MS(\Delta^{IP}_{\alpha}, r) \subset MS(\Delta^{DP}_{\alpha}, \mu).
\eea
\end{proposition}

This proposition indicates that all the functions that can be detected by $\Delta^{IP}_{\alpha}$, can be also be detected by $\Delta^{DP}_{\alpha}$. In other words, the direct test appears to be more efficient in the sense that its associated maxiset is larger. One may ask the question of the strict inclusion. In order to provide an answer, we will consider a specific setting and prove in particular that the inverse testing procedure may miss some functions that can be detected by the direct case.  \\

We are now considering the classical setting of the moderately ill-posed inverse problem, namely we assume that for some $t>0$ and any $k$, $b_k\sim k^{-t}$. We also assume that we are in the case where the calibration is the minimax one. In this case the two terms $r_\varepsilon^2$ and $\sqrt{\sum_{k=1}^{D_\varepsilon} b_k^{-4}}$ are equally balanced, so that we are in the limit case described in Section~\ref{s:3-maxi}.

\begin{proposition}\label{classicalpro}
Let $s,t>0$. Consider the case where for any $\varepsilon \in (0,1)$, $D_\varepsilon \sim \varepsilon^{-\frac{4}{1+4(s+t)}}$. Assume that $r_\varepsilon \sim \varepsilon^{\frac{4s}{1+4(s+t)}}$ and $\mu_\varepsilon \sim b_{D_\varepsilon}r_\varepsilon \sim \varepsilon^{\frac{4(s+t)}{1+4(s+t)}}$ with $b \sim (b_k)_{k\in \BB{N}^*}\sim (k^{-t})_{k\in \BB{N}^*}$.
Then,  there exist functions $\theta$ such that
$$ \theta \in MS^{dec}(\Delta^{DP}_{\alpha}, \mu) \quad \mathrm{but} \quad  \theta \not\in MS^{dec}(\Delta^{IP}_{\alpha}, \mu).$$
\end{proposition}

\begin{remark}
With the choice of operator given in Proposition~\ref{classicalpro}, (\ref{eq:cond}) is clearly satisfied and therefore (\ref{eq:embIPDP}) holds.
\end{remark}

\section{Conclusion}\label{s:5}

In this paper, we adapt the maxiset approach for signal detection in inverse problems. This novel tool for assessing the performance of testing procedures has been exposed to different classical settings. In particular, it allows to compare the so-called direct and inverse approaches. We have established that direct methods are associated to strictly larger maxisets in many cases, which make such testing procedure more interesting for practical purpose. \\

This contribution provides a novel way for researchers to assess the performance of their testing procedures. At the core of our future investigations will be adaptation of our methods to the operator setting leading to the new concept of maxi-class. \\

In order to conclude this discussion, we mention that we are aware of a recent paper of \cite{Ermakov}, that has been recently published in a similar setting while we were finalizing this article. Although this paper also provide a definition of maxiset in testing context, we stress that it uses different constraints on the set of interest. Moreover, it does not consider the inverse problem setting, while we provide a comparison between direct and inverse approach. In our opinion, both contributions are complementary and reveals different aspects of the same problem.

\section{Proofs}\label{s:6}

\subsection{Technical results}
\label{s:tech}

In this section, we recall and slightly extend some results that will be useful in the following. More details regarding these results, e.g., context and extended discussions, can be found in \cite{Baraud} and \cite{LLM2012}.\\

\begin{proposition}\label{prop1}
Fix $(\alpha, \beta) \in (0,1)^{2}$. There exists $C_{\max}(\alpha,\beta)$ and $C_{\min}(\alpha,\beta)$ such that, for all $\varepsilon \in (0,1)$ and $D_\varepsilon \in \BB{N}^*$
\begin{eqnarray*}
(i) & & \sum_{k=1}^{D_\varepsilon} \theta_k^2>C_{\max}(\alpha,\beta) \varepsilon^2 \sqrt{\sum_{k=1}^{D_\varepsilon} b_k^{-4}} \quad \Rightarrow \mathbb{P}_\theta(\Delta_{\alpha,\varepsilon}^{IP} =0)\leq \beta, \\
(ii) & & \sum_{k=1}^{D_\varepsilon} \theta_k^2 \leq  C_{\min}(\alpha,\beta) \varepsilon^2 \sqrt{\sum_{k=1}^{D_\varepsilon} b_k^{-4}} \quad \Rightarrow \mathbb{P}_\theta(\Delta_{\alpha,\varepsilon}^{IP} =0)> \beta.
\end{eqnarray*}
\end{proposition}
\noi{\bf Proof:} We start with the proof of item (i). A more precise proof is provided in \cite{LLM2012}. In particular, the authors take advantage on available results on $\chi^2$ weighted statistics. This allows a better dependency of the constant $C_{\max}(\alpha,\beta)$ w.r.t. the values of $\alpha$ and $\beta$. For the sake of completeness, we reproduce a simpler version of the proof, based on Markov Inequality. Recall that we defined, for any $\varepsilon \in (0,1)$, $T_{D_\varepsilon}$ as $T_{D_\varepsilon}=\sum_{k=1}^{D_\varepsilon}b^{-2}_k(y_k^2-\varepsilon^2)$ and $t_{\alpha,\varepsilon}=C_{\alpha,1}\varepsilon^2\sqrt{\sum_{k=1}^{D_\varepsilon} b_k^{-4}}$ for some constant $C_{\alpha,1} >0$.
Then 
\begin{itemize}
    \item $\mathbb{E}_\theta\left(T_{D_\varepsilon}\right)=\sum_{k=1}^{D_\varepsilon}\theta_k^{2},$
    \item $\mathbb{V}ar_\theta\left(T_{D_\varepsilon}\right)=2\varepsilon^4\sum_{k=1}^{D_\varepsilon}b_k^{-4}+4 \varepsilon^2\sum_{k=1}^{D_\varepsilon}b_k^{-2}\theta_k^{2}.$
\end{itemize} 

 Let $\theta=(\theta_k)_{k\in \mathbb{N}^{*}}$ be such that:
\begin{equation}
    \sum_{k=1}^{D_\varepsilon} \theta_k^2>C_{\max}(\alpha,\beta) \varepsilon^2 \sqrt{\sum_{k=1}^{D_\varepsilon} b_k^{-4}}.
    \label{eq:tail}
\end{equation}

Provided that $C_{\max}(\alpha,\beta)>C_{\alpha,1}$, by using the Bienayme-Chebyshev inequality, one gets
\begin{eqnarray*}
\mathbb{P}_\theta\left(\Delta^{IP}_{\alpha,\varepsilon}=0\right)&=&\mathbb{P}_\theta\left(T_{D_\varepsilon} \leq t_{\alpha, \varepsilon}\right)\\
&=&\mathbb{P}_\theta\left(\mathbb{E}_\theta\left(T_{D_\varepsilon}\right)-T_{D_\varepsilon} \geq \mathbb{E}_\theta\left(T_{D_\varepsilon}\right)-t_{\alpha, \varepsilon}\right)\\
&\leq& \frac{\mathbb{V}ar_\theta\left(T_{D_\varepsilon}\right)}{\left(\mathbb{E}_\theta\left(T_{D_\varepsilon}\right)-t_{\alpha,\varepsilon}\right)^2}\\
&=& \frac{2\varepsilon^4\sum_{k=1}^{D_\varepsilon}b_k^{-4}+4 \varepsilon^2\sum_{k=1}^{D_\varepsilon}b_k^{-2}\theta_k^{2}}{\left(\sum_{k=1}^{D_\varepsilon}\theta_k^{2}-C_{\alpha,1}\varepsilon^2\sqrt{\sum_{k=1}^{D_\varepsilon}b_k^{-4}}\right)^2} \\
& \leq & \frac{2\varepsilon^4\sum_{k=1}^{D_\varepsilon}b_k^{-4}+4 \varepsilon^2 \max_{l=1,...,D_\varepsilon} b_l^{-2} \times \sum_{k=1}^{D_\varepsilon}\theta_k^{2}}{\left(\sum_{k=1}^{D_\varepsilon}\theta_k^{2}-C_{\alpha,1}\varepsilon^2\sqrt{\sum_{k=1}^{D_\varepsilon}b_k^{-4}}\right)^2}\\
& \leq & \frac{2}{(C_{\max}(\alpha,\beta)-C_{\alpha,1})^2}+ \frac{4 \varepsilon^2  \max_{l=1,...,D_\varepsilon} b_l^{-2} C^2_{\max}(\alpha,\beta)}{(C_{\max}(\alpha,\beta)-C_{\alpha,1})^2 \sum_{k=1}^{D_\varepsilon}\theta_k^{2}}\\
& \leq & \frac{2+ 4C_{\max}(\alpha,\beta)}{(C_{\max}(\alpha,\beta)-C_{\alpha,1})^2}\\
&\leq& \beta.
\end{eqnarray*}
\noi The last inequality is obtained if $C_{\max}(\alpha, \beta)$ is large enough. More precisely, we can choose $C_{\max}(\alpha, \beta)$ which depends on both $C_{\alpha,1}$ and $\beta$ as:  \begin{eqnarray}\label{Cmax} \hbox{the maximum of }C_{\alpha,1}\hbox{ and the positive solution of the following equation} \quad \frac{2+ 4C_{\max}(\alpha,\beta)}{(C_{\max}(\alpha,\beta)-C_{\alpha,1})^2}
=\beta.\end{eqnarray}

\begin{remark}
We let the reader be convinced that the smaller $\beta$ the larger the chosen $C_{\min}(\alpha, \beta)$ and $C_{\max}(\alpha, \beta)$.
\end{remark}
%\noi Remark: We use the following assumption: 
%$$\sup_{\varepsilon>0}\frac{\sum_{k=1}^{D_\varepsilon}b_k^{-2}\theta_k^2}{\varepsilon^2\sum_{k=1}^{D_\varepsilon}b_k^{-4}} <\infty.$$

\noindent 
We now prove the item (ii) of Proposition \ref{prop1}. 
%As in the proof of item (i), $t_{\alpha,\varepsilon}$ is chosen such that $t_{\alpha,\varepsilon}=C_{1,\alpha}\varepsilon^2\sqrt{\sum_{k=1}^{D_\varepsilon} b_k^{-4}}$ for some constant $C_{1,\alpha} > 0$. We now give a lower bound of the Type-II error
%\begin{eqnarray*}
%\mathbb{P}_\theta\left(T_{D_\varepsilon}\leq t_{\alpha,\varepsilon}\right)&=&\mathbb{P}_\theta\left(T_{D_\varepsilon}-\mathbb{E}(T_{D_\varepsilon})\leq t_{\alpha,\varepsilon}-\mathbb{E}(T_{D_\varepsilon})\right)\\
%&=&\mathbb{P}_\theta\left(T_{D_\varepsilon}-\mathbb{E}(T_{D_\varepsilon})\leq t_{\alpha,\varepsilon}-\mathbb{E}(T_{D_\varepsilon})\right).
%\end{eqnarray*}
To prove that $\mathbb{P}_\theta\left(\Delta_{\alpha,\varepsilon}^{IP}=0\right)> \beta$ is equivalent to prove that
\begin{equation}
\mathbb{P}_\theta\left(T_{D_\varepsilon}-\mathbb{E}(T_{D_\varepsilon})\geq t_{\alpha,\varepsilon}-\mathbb{E}(T_{D_\varepsilon})\right)< 1-\beta.
\label{eq:upper-bound}
\end{equation}
To show this inequality, we apply Lemma 2 of Laurent et al. (2012) with $x:=-\log(1-\beta)$ and $\sigma_k \equiv \varepsilon b_k^{-1}$. Setting $\Sigma=\varepsilon^4\sum_{k=1}^{D_\varepsilon} b_k^{-4}+2\varepsilon^2 \sum_{k=1}^{D_\varepsilon} b_k^{-2}\theta_k^2$, we then get that~\eqref{eq:upper-bound} holds provided that 
\[
t_{\alpha,\varepsilon}-\mathbb{E}(T_{D_\varepsilon})\geq 2\sqrt{\Sigma x}+2\max_{l=1,...,D_\varepsilon} b_l^{-2}\varepsilon^2 x.
\]
Observe that
\begin{eqnarray*}
\sqrt{\Sigma}&\leq& \varepsilon^2\sqrt{\sum_{k=1}^{D_\varepsilon} b_k^{-4}}+\sqrt{2\varepsilon^2 \sum_{k=1}^{D_\varepsilon} b_k^{-2}\theta_k^2}\\
&\leq& \varepsilon^2\sqrt{\sum_{k=1}^{D_\varepsilon} b_k^{-4}}+\sqrt{2}\varepsilon \max_{l=1,...,D_{\varepsilon}} b_l^{-1}\sqrt{\sum_{k=1}^{D_\varepsilon} \theta_k^2}\\
&\leq& \left(1+\sqrt{2C_{\min}(\alpha,\beta)}\right)\varepsilon^2\sqrt{\sum_{k=1}^{D_\varepsilon} b_k^{-4}}\;,
\end{eqnarray*}
the last inequality comes from the assumption
\[
\sum_{k=1}^{D_\varepsilon} \theta_k^2\leq C_{\min}(\alpha,\beta) \varepsilon^2\sqrt{\sum_{k=1}^{D_\varepsilon} b_k^{-4}}.
\]
Since $\beta <1/2$,  $x<1$ (and therefore $x < \sqrt{x}$). So, a sufficient condition to get the inequality $t_{\alpha,\varepsilon}-\mathbb{E}(T_{D_\varepsilon})\geq 2\sqrt{\Sigma x}+2\max_{l=1,...,D_\varepsilon} b_l^{-2}\varepsilon^2 x$ is 
\[
C_{\alpha, 1} \varepsilon^2\sqrt{\sum_{k=1}^{D_\varepsilon} b_k^{-4}}-\sum_{k=1}^{D_\varepsilon} \theta_k^2\geq \left(2\left(1+\sqrt{2C_{\min}(\alpha,\beta)}\right)\sqrt{\sum_{k=1}^{D_\varepsilon}b_k^{-4}}\right) \varepsilon^2\sqrt{-\log(1-\beta)} + 2b_{D_{\varepsilon}}^{-2}\varepsilon^2 \sqrt{-\log(1-\beta)}.
\]
This relationship is satisfied if $C_{\min}(\alpha,\beta)$ is small enough. The choice \begin{eqnarray}\label{Cmin}
C_{\min}(\alpha,\beta)= \left(\sqrt{-2\log(1-\beta)+(C_{\alpha, 1}-4\sqrt{-\log(1-\beta)})}-\sqrt{-2\log(1-\beta)}\right)^{\frac{1}{2}}\end{eqnarray} can be done in that context. This finishes the proof since this condition is compatible with $\Delta^{IP}_{\alpha,\varepsilon}$ is a $\alpha$ level test for a chosen $C_{\alpha, 1}$ sufficiently large enough.
\begin{flushright}
$\Box$
\end{flushright}

\begin{remark}
The smaller $\beta$ the more difficult the signal detection problem. This difficulty of the signal detection problem is highlighted here in the value of $C_{\min}(\alpha, \beta)$. Indeed, we let the reader be convinced that the smaller $\beta$ the larger the chosen $C_{\min}(\alpha, \beta)$ and therefore the larger the $l_2$ norm of the  the first terms of the signal $\theta$ have to be to ensure that $\mathbb{P}_\theta(\Delta_{\alpha,\varepsilon}^{IP} =0)\leq \beta$ where $\varepsilon \in (0,1)$.
\end{remark}

In the same spirit of Proposition \ref{prop1}, the following result holds too.

\begin{proposition}\label{prop2}
There exists $C'_{\min}(\alpha,\beta)$ and $C'_{\max}(\alpha,\beta)$ such that, for all $\varepsilon \in (0,1)$
\begin{eqnarray*}
(i) & & \sum_{k=1}^{D_\varepsilon} b_k^2\theta_k^2>C'_{\max}(\alpha,\beta) \varepsilon^2 \sqrt{D_\varepsilon} \quad \Rightarrow \mathbb{P}_\theta(\Delta_{\alpha,\varepsilon}^{DP} =0)\leq \beta, \\
(ii) & & \sum_{k=1}^{D_\varepsilon} b_k^2\theta_k^2< C'_{\min}(\alpha,\beta) \varepsilon^2 \sqrt{D_\varepsilon} \quad \Rightarrow \mathbb{P}_\theta(\Delta_{\alpha,\varepsilon}^{DP} =0)> \beta.
\end{eqnarray*}
\end{proposition}
\noi Since the proof of Proposition \ref{prop2} is analogous to the one of Proposition \ref{prop1}, we omit it. We precise that the values of $C'_{\min}(\alpha,\beta)$ and $C'_{\max}(\alpha,\beta)$ are obtained from the ones of $C_{\min}(\alpha,\beta)$ and $C_{\max}(\alpha,\beta)$ by replacing $C_{\alpha,1}$ by $C_{\alpha,2}$.
\begin{flushright}
$\Box$
\end{flushright}
\vspace{-0.5cm}
\subsection{Proof of Theorem \ref{maintheo1}}\label{s:proof-maintheo1}
We first prove the following embedding property : $\mathcal{F}_{r,D}(C_{\max}(\alpha,\beta))\subset MS(\Delta_{\alpha}^{IP},r)$ and thereafter we prove that the embedding property $MS(\Delta_{\alpha}^{IP}, r)\subset \mathcal{F}_{r,D}(C_{\min}(\alpha,\beta))$, whatever the choice of the detection rate $r=(r_\varepsilon)_{\varepsilon>0}$.\\

First we prove that  $\mathcal{F}_{D}(C_{\max}(\alpha,\beta)) \subset MS(\Delta_{\alpha}^{IP}, r)$. Fix $\varepsilon \in (0,1)$ and let $\theta \in  \mathcal{F}_{D}(C_{\max}(\alpha,\beta))$ satisfying $\|\theta\|^2\geq r_\varepsilon^2$. Then
\[
 \sum_{k=1}^{D_\varepsilon}\theta_k^2>C_{\max}(\alpha,\beta)\sqrt{\sum_{k=1}^{D_\varepsilon}b_k^{-4}}\Rightarrow\mathbb{P}_\theta(\Delta_{\alpha, \varepsilon}=0)\leq \beta.
\]
This assertion is obvious thanks to item (i) of Proposition \ref{prop1}. \\

Conversely, we shall prove that $MS(\Delta_{\alpha}^{IP}, r)\subset \mathcal{F}_{D}(C_{\min}(\alpha,\beta))$, or equivalently that 
$$ \theta \not \in \mathcal{F}_{D}(C_{\min}(\alpha,\beta)) \Rightarrow \theta \not \in MS(\Delta_{\alpha}^{IP},r).$$
Assume that $\theta \not \in \mathcal{F}_{D}(C_{\min}(\alpha,\beta)) $. Then there exists $\varepsilon_0 \in (0,1)$ such that 
$$
\|\theta\|^2\geq r_{\varepsilon_0}^2\mbox{ and }\sum_{k=1}^{D_{\varepsilon_0}}\theta_k^2 \leq C_{\min}(\alpha,\beta)\varepsilon_0^2\sqrt{\sum_{k=1}^{D_{\varepsilon_0}}b_k^{-4}}.
$$
Following item $(ii)$ of Proposition \ref{prop1}, one deduces that
\begin{equation}\label{eq:non-contr}
\|\theta\|^2\geq r_{\varepsilon_0}^2\mbox{ and } \mathbb{P}_\theta(\Delta_{\alpha, \varepsilon_0}=0) > \beta.
\end{equation}
Therefore we immediately deduce that $\theta \not \in MS(\Delta_{\alpha}^{IP},r)$. \\

The proof of the second part of Theorem \ref{maintheo1} follows exactly the same lines and is left to the interested reader.

\subsection{Proof of Theorem \ref{maintheo2}}\label{s:proof-maintheo2}
As in the previous proof, we concentrate our attention on the maxiset associated to the sequence of testing procedures $(\Delta_{\alpha,\varepsilon}^{IP})_{\varepsilon>0}$. We recall that our aim is to prove that the robust maxiset $MS^{dec}(\Delta_{\alpha}^{IP},r)$ can be identified, up to a constant, to the set $\CAL{F}^{dec}_{r, D}(C)$ defined as  
\beas
\CAL{F}^{dec}_{r, D}(C) = \left\{\theta\in l_2(\mathbb{N}^{*}), \ \forall \varepsilon \in (0,1) ; \sum_{k>D_\varepsilon}\theta_k^2 \leq r_{\varepsilon}^2 - C\varepsilon^2\sqrt{\sum_{1}^{D_\varepsilon}b_k^{-4}} \right\},
\eeas
We begin the proof with the first inclusion that is $\CAL{F}^{dec}_{r, D}(C_{\max}(\alpha,\beta)) \subset MS^{dec}(\Delta_{\alpha}^{IP},r)$. It can be easily seen that for all $\varepsilon \in (0,1)$ 
\begin{eqnarray*}
\| \theta \|^2 \geq r_\varepsilon^2 \; \mathrm{and} \; \theta \in  \CAL{F}^{dec}_{r, D}(C_{\max}(\alpha,\beta))
& \Rightarrow & \sum_{k> D_\varepsilon} \theta_k^2 < \| \theta \|^2 - C_{\max}(\alpha,\beta)\varepsilon^2\sqrt{\sum_{1}^{D_\varepsilon}b_k^{-4}} \\
& \Rightarrow & \sum_{k=1}^{D_\varepsilon} \theta_k^2 > C_{\max}(\alpha,\beta)\varepsilon^2\sqrt{\sum_{1}^{D_\varepsilon}b_k^{-4}} \\
& \Rightarrow & \mathbb{P}_\theta(\Delta^{IP}_{\varepsilon,\alpha}=0) \leq \beta \: \text{ because of item (i) of Proposition }\ref{prop1}.
\end{eqnarray*}
Since, one can easily check that $\CAL{F}^{dec}_{r, D}(C_{\max}(\alpha,\beta))$ is stable by decimation, this entails that $\CAL{F}^{dec}_{r, D}(C_{\max}(\alpha,\beta)) \subset MS^{dec}(\Delta_{\alpha}^{IP},r)$.\\ %We stress that the decimation constraint is not necessary in order to get this inclusion. \\

Now, we turn our attention to the proof of the second inclusion that is $MS^{dec}(\Delta_{\alpha}^{IP},r) \subset \CAL{F}^{dec}_{\sqrt{2}r, D}(C_{\min}(\alpha,\beta))$. We consider the situation where for any $\varepsilon \in (0,1)$
$$ C_{\max}(\alpha,\beta)\varepsilon^2 \sqrt{\sum_{k=1}^{D_\varepsilon} b_k^{-4}} < r_{\varepsilon}^2$$ 
(otherwise $\CAL{F}^{dec}_{r, D}(C_{\max}(\alpha,\beta))$ is empty)
% and assume, without loss of generality that 
%$$C_{\min}(\alpha,\beta)\varepsilon^2 \sqrt{\sum_{k=1}^{D_\varepsilon} b_k^{-4}} < r_{\varepsilon}^2/2.$$
Let $\theta \not \in \CAL{F}^{dec}_{\sqrt{2}r, D}(C_{\min}(\alpha,\beta))$. In particular, since clearly $C_{\max}(\alpha,\beta)\geq C_{\min}(\alpha,\beta)$, there exists $\varepsilon_1 \in (0,1)$ such that 
$$ \sum_{k>D_{\varepsilon_1}} \theta_k^2 \geq 2r_{\varepsilon_1}^2 - C_{\min}(\alpha,\beta) \varepsilon_1^2 \sqrt{\sum_{k=1}^{D_{\varepsilon_1}} b_k^{-4}} \geq 2r_{\varepsilon_1}^2 - C_{\max}(\alpha,\beta) \varepsilon_1^2 \sqrt{\sum_{k=1}^{D_{\varepsilon_1}} b_k^{-4}} > r_{\varepsilon_1}^2.$$
Hence, the sequence $(\theta_k^{(-D_{\varepsilon_1})})_{k\in\mathbb{N}^{*}}$ satisfies both
$$ \| \theta^{(-D_{\varepsilon_1})}\|^2 = \sum_{k>D_{\varepsilon_1}} \theta_k^2 > r_{\varepsilon_1}^2 \quad \mathrm{and} \quad \mathbb{P}_{\theta^{(-D_{\varepsilon_1})} }(\Delta_{\alpha, \varepsilon_1}^{IP} = 0) > \beta,$$
where the last inequality is obtained thanks to item (ii) of Proposition \ref{prop1} and the fact that $\theta^{(-D_{\varepsilon_1})}_k=0$ for all $k\leq D_{\varepsilon_1}$. This entails that $\theta^{(-D_{\varepsilon_1})} \not \in MS(\Delta_{\alpha}^{IP}, r)$ and therefore $\theta\not \in MS^{dec}(\Delta_{\alpha}^{IP}, r)$ because of the decimation constraint. Hence, we get 
$$ \theta \not \in \CAL{F}^{dec}_{\sqrt{2}r, D}(C_{\min}(\alpha,\beta)) \Rightarrow \theta\not \in MS(\Delta_{\alpha}^{IP},r)$$
which implies that
$$ MS^{dec}(\Delta_{\alpha}^{IP},r) \subset  \CAL{F}^{dec}_{\sqrt{2}r, D}(C_{\min}(\alpha,\beta)) .$$
This concludes the proof. 
\begin{flushright}
$\Box$
\end{flushright}

%\subsection{Proof of Proposition~\ref{pro:emb}}\label{s:proof_emb}
%\noindent
%We consider the case where $\Delta^{(.)}_{\alpha}$ is $\Delta_{\alpha}^{IP}$, the other case being similar. In view of Theorem~\ref{maintheo2}, it is sufficient to prove that for any $C,C'>0$, $\,\mathcal{F}_{r,D}^{dec}(C)\subset\mathcal{F}_{r',D}^{dec}(C')$. And to do so we only have to check that for $\varepsilon$ small enough
%\[
%r^{'2}_\varepsilon-r_\varepsilon^2-C' \varepsilon^2\sqrt{\sum_{k=1}^{D_\varepsilon}b_k^{-4}}\geq 0.
%\]
%Since $r_\varepsilon=o(r'_\varepsilon)$, it is sufficient to prove that for $\varepsilon$ small enough
%\[
%r^{'2}_\varepsilon-C' \varepsilon^2\sqrt{\sum_{k=1}^{D_\varepsilon}b_k^{-4}}> 0.
%\]
%We remark that if it was wrong it would mean that both $\mathcal{F}_{r,D}^{dec}(C)$ and $\mathcal{F}_{r',D}^{dec}(C')$ are empty sets which implies also the result.
%\begin{flushright}
%$\Box$
%\end{flushright}

\subsection{Proof of Proposition~\ref{maincor}}\label{s:proof-maincor}
The proof of Proposition~\ref{maincor} directly follows Theorem~\ref{maintheo2} and from the fact that 
\[
\CAL{F}_{r,D}^{dec}(C_{\min}(\alpha,\beta))\subset\CAL{G}_{\mu,D}^{dec}(C'_{\max}(\alpha,\beta)).
\]
Let us prove this embedding of sequence spaces. Assume first that $\theta\in \CAL{F}_{r,D}^{dec}(C_{\min}(\alpha,\beta))$, then according to assumption (\ref{eq:cond}) on the sequence $(b_k)_{k\in \mathbb{N}^*}$, for any $\varepsilon \in (0,1),$
\[
\sum_{k>D_\varepsilon}b_k^2 \theta_k^2 \leq b_{D_\varepsilon}^2\left(\sum_{k>D_\varepsilon}\theta_k^2 \right)< b_{D_\varepsilon}^2 \left(r_\varepsilon^2-C_{\min}(\alpha,\beta)\varepsilon^2\sqrt{\sum_{k=1}^{D_\varepsilon}b_k^{-4}}\right).
\]
Since $\mu_\varepsilon=b_{D_\varepsilon}r_\varepsilon$ and $b_{D_\varepsilon}^2\sqrt{\sum_{k=1}^{D_\varepsilon}b_k^{-4}}=\sqrt{\sum_{k=1}^{D_\varepsilon}(b_{D_\varepsilon}/b_k)^{4}}\geq \frac{C'_{\max}(\alpha,\beta)}{C_{\min}(\alpha,\beta)}\sqrt{D_\varepsilon}$, we get
\[
\sum_{k>D_\varepsilon}b_k^2 \theta_k^2 <  \mu_\varepsilon^2-C'_{\max}(\alpha,\beta)\varepsilon^2\sqrt{D_\varepsilon}
\]
which exactly means that $\theta\in \CAL{G}_{\mu,D}^{dec}(C'_{\max}(\alpha,\beta))$.

\subsection{Proof of Proposition~\ref{classicalpro}}\label{s:proof-classicalpro}
We consider the case where for any $k\in\BB{N}^*$, $b_k\sim k^{-t}$. One has 
\[
r_\varepsilon^2\sim\varepsilon^2\sqrt{\sum_{k=1}^{D_\varepsilon}b_k^{-4}}\sim \varepsilon^2 D_\varepsilon^{(1+4t)/2}\sim D_\varepsilon^{-(1+4(s+t))/2}D_\varepsilon^{(1+4t)/2}=D_\varepsilon^{-2s}.
\]
We are then in the limit case. According to the value of $C$, $\mathcal{F}_{r,D}(C)$ may be empty or may be as follows:
\[
\mathcal{F}^{dec}_{r,D}(C):=\left\{\theta: \ \forall \varepsilon \in (0,1), \ \sum_{k>D_\varepsilon}\theta_k^2\leq CD_\varepsilon^{-2s}\right\}=\left\{\theta:\sup_{K\in \mathbb{N}^*} K^{2s}\sum_{k>K}\theta_k^2\leq C\right\}
\]
which is a Besov space. Note that since the maxiset is only between the two sets $\CAL{F}_{r,D}^{dec}(C_1)$ and $\CAL{F}_{r,D}^{dec}(C_2)$ with $C_1$ and $C_2$ unknown constants we cannot give an equality concerning the maxiset. Using similar computations, we get that 
$$ \mathcal{G}_{\mu,D}^{dec}(C) = \left\lbrace \theta:\sup_{K\in \mathbb{N}^*} K^{2(s+t)}\sum_{k>K}b_k^2 \theta_k^2\leq C\right\}.$$
Our aim below is to exhibit a sequence $\theta$ such that $\theta \in  \mathcal{G}^{dec}_{\mu,D}(C) $ but $\theta \not \in \mathcal{F}^{dec}_{\mu,D}(C')$, whatever the value of the constant $C'$. To this end, let us consider the sequence $\theta=(\theta_k)_{k\in \mathbb{N}^*}$ such that 
$$ \theta_k = \frac{2^{-js}}{k} \quad \forall k\in \lbrace 2^j,\dots, 2^{j+1}-1 \rbrace, \ j\in \mathbb{N}^*. $$
Let $K\in \mathbb{N}^*$ be fixed, and $j_0$ such that $ 2^{j_0} \leq K \leq 2^{j_0+1}$. Then, we can check that 
$$ \sum_{k>K} b_k^2\theta_k^2 \lesssim\sum_{k>K} k^{-(2t+2)} 2^{-2js} \lesssim 2^{-2j_0s} K^{-2t} \sim K^{-2(s+t)}.$$
In the same time, if $K=2^{j_0}$, 
$$ \sum_{k>K} \theta_k^2 \geq \sum_{k=2^{j_{0}}}^{2^{j_0+1}-1} \frac{2^{-2js}}{k^2} \sim 2^{-2j_0s} \ln(2^{j_0+1}-2^{j_{0}}) \sim 2^{-2j_0s} j_0 >> K^{-2s},$$
when $j$ increases.

%\subsection{Proof of the results of Section~\ref{s:4}}

%\section{Perspectives}
%\subsection{Anisotropic setting}
%We assume that we are given an anisotropic function f and an operator A whose anistropy may
%differ from that of f . We want to tune the parameter
%\subsection{Unknown operator}
%\begin{itemize} 
%\item Noise on eigenvalues
%\item Family of operators: agregation of operators (Bonferoni)
%\end{itemize}

%\biblist

%\section*{Acknowledgments}
\end{document}